\documentclass[12pt, dvips]{amsart}

\swapnumbers\newtheorem{theorem}{Theorem}[section]

\newtheorem{proposition}[theorem]{Proposition}

\newtheorem{lemma}[theorem]{Lemma}

\newtheorem{definition}[theorem]{Definition}

\newtheorem{remark}[theorem]{Remark}

\newtheorem{example}[theorem]{Example}

\numberwithin{equation}{section}

\usepackage{amssymb,latexsym, amsmath, amscd, bbm, graphicx, epsfig, psfrag}

\def\s{\smallskip}
\def\m{\medskip}
\def\ni{\noindent}
\def\ra{\rightarrow}
\def\ha{\hookrightarrow}

\def\RR{\mathbbm{R}}
\def\NN{\mathbbm{N}}

\def\gg{\gamma}

\def\ee{\epsilon}

\def\ff{\varphi}
\def\ll{\lambda}
\def\oo{\omega}

\def\gt{\vartheta}

\def\ca{{\mathcal A}}

\def\cl{{\mathcal L}}

\def\cv{{\mathcal V}}

\def\supp{\operatorname{supp}\:\!}
\def\length{\operatorname{length}\:\!}

\newcommand{\proofend}{\hspace*{\fill} $\Box$\\}
\newcommand{\diam}{\hspace*{\fill} $\Diamond$}

\begin{document}

\title[]{An extension theorem in symplectic geometry}

\date{\today}

\author{Felix Schlenk}
\address{(F.\ Schlenk) ETH Z\"urich, CH-8092 Z\"urich, Switzerland}
\email{felix@math.ethz.ch}

\begin{abstract}  
We extend the ``Extension after Restriction Principle'' for symplectic
embeddings of bounded starlike domains to a large class of symplectic
embeddings of unbounded starlike domains.
\end{abstract}

\maketitle

\markboth{{\rm An extension theorem in symplectic geometry}}{{}}

\section{Introduction}  
\ni
We endow each open subset $U$ of Euclidean space $\RR^{2n}$ 
with the standard symplectic form
\[
\oo_0 \,=\, \sum_{i=1}^n dx_i \wedge dy_i .
\]
A smooth embedding $\ff \colon U \ha \RR^{2n}$ is called symplectic if 
$\ff^* \oo_0 \,=\, \oo_0$.
In particular, every symplectic embedding preserves the volume form
$\tfrac{1}{n!} \omega_0^n$ and hence the Lebesgue measure on $\RR^{2n}$.
Recall that a domain in $\RR^{2n}$ is by definition a non-empty open
connected subset of $\RR^{2n}$.

\begin{definition}
{\rm
Consider a symplectic embedding $\ff \colon U \ha \RR^{2n}$ 
of a domain $U$ in $\RR^{2n}$.
We say that the pair $(U, \ff)$ has the {\it extension property}\,
if for each subset $A \subset U$ whose closure in $\RR^{2n}$ is
contained in $U$ there exists a symplectomorphism $\Phi_A$ of $\RR^{2n}$ 
such that $\Phi_A |_A = \ff |_A$.
}
\end{definition}

Not every pair $(U, \ff)$ as above has the extension property as the following
example shows.

\begin{example}  \label{exa:0}
{\rm
For $0 \le r_0<r_1<\infty$ we define the open annulus 
\[
A \left( r_0,r_1 \right) \,=\,
\left\{ (x,y) \in \RR^2 \mid r_0^2 < x^2+y^2 < r_1^2 \right\} .
\]
Let $\ff \colon A(0,3) \ra A(4,5) \subset \RR^2$ be the symplectic embedding which
in polar coordinates is given by
\[
\ff (r, \gt) \,=\, \left( \sqrt{r^2+16}, \gt \right) .
\]
Any smooth extension of $\ff |_{A(1,2)}$ to $\RR^2$ maps the disc of
area $\pi$ to the disc of area $17\pi$ and hence cannot be symplectic.
\diam
}
\end{example}

On the other hand, the well-known 
``Extension after Restriction Principle'' \cite{EH1}, 
which is reproved below, states that a pair $(U, \ff)$ has the
extension property whenever the geometry of $U$ is simple enough.
Recall that a subset $U$ of $\RR^{2n}$ is said to be {\it starlike}\;\! if 
$U$ contains a point $p$ such that for every point $x \in U$ the
straight line between $p$ and $x$ is contained in $U$.
 
\begin{proposition}\  \label{pa2:ear}
{\rm (Extension after Restriction Principle)}
Assume that $\ff \colon U \ha \RR^{2n}$ is a symplectic embedding of a
bounded starlike domain $U \subset \RR^{2n}$.
Then the pair $(U,\ff)$ has the extension property. In fact, for any subset $A \subset U$ whose closure in $\RR^{2n}$ is
contained in $U$ there exists a compactly 
supported symplectomorphism $\Phi_A$ of $\RR^{2n}$ such that
$\Phi_A |_A = \ff |_A$.
\end{proposition}

The purpose of this paper is to prove the extension property for a
large class of symplectic embeddings of unbounded starlike domains.
The following example shows that it is not enough to assume that $U$
is starlike.

\begin{example}  \label{exa:1}
{\rm
We let $U \subset \RR^2$ be the strip 
$]1, \infty[ \,\times\, ]\:\!\!-1,1[$.
Combining the methods used in Step 1 and Step 4 of \text{Section 2.2}
in \cite{Diss} 
we find a symplectic embedding $\ff \colon U \ha \RR^2$ such that
$\ff (k,0) = \left( \tfrac{1}{k}, 0 \right)$, $k = 2,3, \dots$.
Then there does not exist any subset $A$ of $U$ containing the set 
$\{ (k,0) \mid k = 2,3, \dots \}$
for which $\ff |_A$ extends to a diffeomorphism of $\RR^2$.
\diam
}
\end{example}

Observe that if $(U, \ff)$ has the extension property, then $\ff$ is
proper in the sense that each subset $A \subset U$
whose closure in $\RR^{2n}$ is contained in $U$ and whose image $\ff
(A)$ is bounded is bounded. 
The map $\ff$ in Example \ref{exa:1} is not proper in this sense. 
However, the map $\ff$ in the following example is proper in this sense, 
and still $(U, \ff)$ does not have the extension property.

\begin{example}  \label{exa:2}
{\rm
Let $U \subset \RR^2$ be the strip $\RR \,\times\, ]\;\!\!-1,0\:\![$,
and let 
\[
A \,=\, \left\{ (x,y) \in U \;\big|\:\! \left| y+\tfrac{1}{2} \right| \le f(x)
\right\}
\]
where $f \colon \RR \ra \left] 0,\tfrac{1}{2} \right[$ is a smooth
function such that
\begin{equation}  \label{ea:intf}
\int_\RR \left( \tfrac{1}{2} - f(x) \right) dx \,<\, \infty ,
\end{equation}
\begin{figure}[h] 
 \begin{center}
  \psfrag{0}{$0$}
  \psfrag{1}{$1$}
  \psfrag{-1}{$-1$}
  \psfrag{x}{$x$}
  \psfrag{y}{$y$}
  \psfrag{U}{$U$}
  \psfrag{A}{$A$}
  \psfrag{ff}{$\ff$}
  \leavevmode\epsfbox{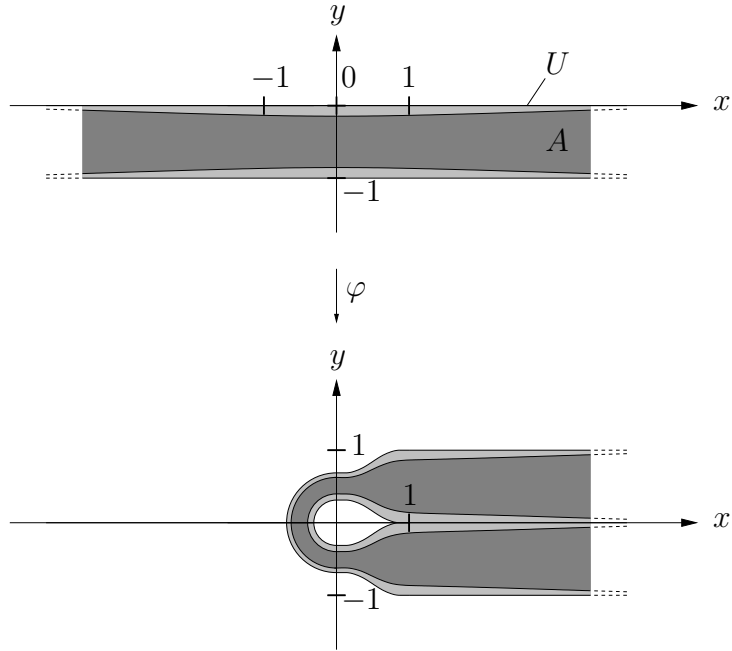}
 \end{center}
 \caption{A pair $(U, \ff)$ which does not have the extension property.}
 \label{figurea1.fig}
\end{figure}
%
%magnification: 38
%

\noindent
cf.\ Figure \ref{figurea1.fig}.
Using the method used in Step 4 of \text{Section 2.2} in \cite{Diss}
we find a symplectic embedding $\ff \colon U \ha \RR^2$ such that
\[
\ff (x,y) = (x,y) \,\text{ if } x \ge 1 
\quad \text{ and } \quad 
\ff (x,y) = (-x,-y) \,\text{ if } x \le -1 ,
\]
cf.\ \text{Figure \ref{figurea1.fig}}.
In view of the estimate \eqref{ea:intf} the component $C$ of $\RR^2
\setminus \ff(A)$ which contains the point $(1,0)$ has finite volume. 
Any symplectomorphism $\Phi_A$ of $\RR^2$ such that $\Phi_A |_A = \ff
|_A$ would map the ``upper'' component of $\RR^2 \setminus A$,
which has infinite volume, to $C$.
This is impossible.
\diam
}
\end{example}

Example \ref{exa:2} shows that the assumption \eqref{ea2:2} on $\ff$ 
in Theorem \ref{ta2:ear} below cannot be omitted.
For technical reasons in the proof of Theorem \ref{ta2:ear} we shall
also impose a mild convexity condition on the starlike domain $U$.
The length of a smooth curve $\gg \colon [0,1] \ra \RR^n$ is defined
by
\[
\length (\gg) \,:=\, \int_0^1 \left| \gg' (s) \right| ds .
\]
On any domain $U \subset \RR^n$ we define a distance function $d_U
\colon U \times U \ra \RR$ by
\[
d_U (z,z') \,:=\, \inf \left\{ \length (\gg) \right\} 
\]
where the infimum is taken over all smooth curves $\gg \colon [0,1] \ra
U$ with $\gg(0) = z$ and $\gg(1) = z'$.
Then $\left| z-z' \right| \le d_U (z,z')$ for all $z,z' \in U$.

\begin{definition}  \label{defa:lip}
{
\rm
We say that a domain $U \subset \RR^{n}$ is a {\it Lipschitz domain}\,
if there exists a constant $\ll >0$ such that
\[
d_U (z,z') \,\le\, \ll \left| z-z' \right|
                        \quad \text{for all }\, z, z' \in U.
\]
}
\end{definition}

Each {\it convex}\, domain $U \subset \RR^{n}$ is a Lipschitz domain
with Lipschitz constant $\ll=1$.
It is not hard to see that there do exist starlike domains
which are not Lipschitz domains.
But we do not know of a starlike domain with smooth boundary which is 
not a Lipschitz domain.

\begin{theorem}  \label{ta2:ear}
Assume that $\ff \colon U \hookrightarrow \RR^{2n}$ is a symplectic
embedding of a
starlike Lipschitz domain $U \subset \RR^{2n}$ such that
there exists a constant $L>0$ satisfying
\begin{equation}  \label{ea2:2}
\left| \ff (z) - \ff (z') \right| \,\ge\, L \left| z-z' \right| \quad \text{for all }\, z, z' \in U.
\end{equation}  
Then the pair $(U,\ff)$ has the extension property.
\end{theorem}

Theorem \ref{ta2:ear} is applied in \cite{Diss} 
to extend a symplectic vanishing theorem for bounded domains to
certain unbounded domains.

\subsection*{Acknowledgements}
I cordially thank Urs Lang, Fran\c cois Laudenbach and Edi Zehnder
for useful discussions.

\section{Proofs}  

\ni
We shall first proceed along the lines of \cite{EH1} and then verify
that our assumptions on $U$ and $\ff$ are
sufficient to push the arguments through.

\m
\noindent
{\bf Step 1. Reduction to a simpler case}

\s
\noindent
We start with observing that we may assume that $U$ is starlike with
respect to the origin and that $\ff (0) = 0$ and $d \ff (0) = id$.
Indeed, suppose that \text{Theorem \ref{ta2:ear}} holds in this
situation, that $A$ is a subset of $U$ whose closure in $\RR^{2n}$ is
contained in $U$, 
and that $U$ is starlike with respect to $p \neq 0$ or that
$\ff (p) \neq 0$ or that $D := d \ff (p) \neq id$.
For $w \in \RR^{2n}$ we denote by $\tau_w$
the translation $z \mapsto z + w$. 
We define the symplectic embedding 
$\psi \colon \left( D \circ \tau_{-p} \right)(U) \ha \RR^{2n}$ by
\[
\psi \,:=\, \tau_{-\ff(p)} \circ \ff \circ \tau_p \circ D^{-1} .
\]
Then $\psi (p) =0$ and $d \psi (p) = id$.
Since $U$ is starlike with respect to $p$ and $D$ is linear, 
the domain $\left( D \circ \tau_{-p} \right)(U)$ is starlike with
respect to the origin, and $\left( D \circ \tau_{-p} \right)(A)$ 
is a subset of $\left( D \circ \tau_{-p} \right)(U)$ whose closure in
$\RR^{2n}$ is contained in $\left( D \circ \tau_{-p} \right)(U)$. 
Assume next that $U$ is a $\ll$-Lipschitz domain.
We fix $w,w' \in \left( D \circ \tau_{-p} \right)(U)$ and set
$z = D^{-1}(w)+p$, $z' = D^{-1}(w')+p$. 
Given any smooth path $\gg \colon [0,1] \ra U$ with $\gg(0) =z$, $\gg(1) 
= z'$, the smooth path 
\[
D \circ \tau_{-p} \circ \gg \colon [0,1] \,\ra\, \left( D \circ \tau_{-p}
\right)(U)
\]
runs from $w$ to $w'$,
and so
\[
d_{\left( D \circ \tau_{-p} \right)(U)}(w,w') \,\le\, \int_0^1 \left| D\;\!
\gg'(s) \right| ds 
\,\le\, \left\| D \right\| \int_0^1 \left| \gg'(s) \right| ds.
\]
It follows that
\begin{eqnarray*}
d_{\left( D \circ \tau_{-p} \right)(U)}(w,w') 
  &\le& \left\| D \right\| d_U (z,z') \\
  &\le& \left\| D \right\| \ll \left| z-z' \right| \\
  &\le& \left\| D \right\| \ll \left\| D^{-1} \right\| \left| w-w' \right|. 
\end{eqnarray*}
Since $w, w' \in \left( D \circ \tau_{-p} \right)(U)$ were arbitrary, we 
conclude that $\left( D \circ \tau_{-p} \right)(U)$ is a
$\left\| D \right\| \left\| D^{-1} \right\| \ll$-Lipschitz domain.
Finally, the assumption \eqref{ea2:2} on $\ff$ yields
\[
\left| \psi (z) - \psi (z') \right| \,\ge\, 
L \left| D^{-1}(z-z') \right| \,\ge\, \frac{L}{\| D \|} |z-z'|
\]
for all $z,z' \in \left( D \circ \tau_{-p} \right)(U)$.
By assumption we therefore find a symplectomorphism 
$\Psi_{\left( D \circ \tau_{-p} \right)(A)}$ of $\RR^{2n}$
such that $\Psi_{\left( D \circ \tau_{-p} \right)(A)} |_{\left( D \circ \tau_{-p} \right)(A)} = \psi |_{\left( D \circ \tau_{-p} \right)(A)}$.
Define the symplectomorphism $\Phi_A$ of $\RR^{2n}$ by
\[
\Phi_A \,:=\, \tau_{\ff(p)} \circ \Psi_{\left( D \circ \tau_{-p}
                                   \right)(A)} \circ D \circ \tau_{-p} . 
\]
Then $\Phi |_A = \ff |_A$, as required.

\m
\noindent
{\bf Step 2. The classical approach}

\s
\noindent
So assume that $U$ is starlike with respect to the origin and that $\ff
(0) = 0$ and $d \ff (0) = id$.
We denote the set of symplectic embeddings of $U$ into $\RR^{2n}$ by 
$\mbox{Symp}\, (U, \RR^{2n})$.
Since $U$ is starlike with respect to the origin we can define a continuous
path $\ff_t \subset \mbox{Symp}\, (U, \RR^{2n})$ by setting
\begin{equation}  \label{de:phi}
\ff_t(z) \,:=\, 
\left\{ \begin{array}{ll}
          z                       & \text{ if }\; t =0 , \\ [0.1em]
          \frac{1}{t} \, \ff (tz) & \text{ if }\; t \in\;]0,1] .
        \end{array}
   \right. 
\end{equation}
The path $\ff_t$ is smooth except possibly at $t=0$. 
In order to smoothen $\ff_t$, we define the diffeomorphism $\eta$ of $[0,1]$
by
\begin{equation}  \label{de:eta}
\eta (t) \,:=\, 
\left\{ \begin{array}{ll}
          0            & \text{ if }\; t = 0 , \\ [0.1em]
          e^2 \;\!e^{-2/t} & \text{ if }\; t \in\;]0,1] ,
        \end{array}
   \right. 
\end{equation}
where $e$ denotes the Euler number, 
and for $t \in [0,1]$ and $z \in U$ we set
\begin{equation}  \label{de:phie}
\phi_t(z) \,:=\, \ff_{\eta(t)} (z) .
\end{equation}
Then $\phi_t$ is a smooth path in $\mbox{Symp}\, (U, \RR^{2n})$. 
We have $\phi_0 = id_U$ and $\phi_1 = \ff$. 

Since $U$ is starlike, it is contractible, and so the same holds true for
all the open sets $\phi_t (U)$, $t \in [0,1]$.
We therefore find a smooth time-dependent Hamiltonian function 
\begin{equation}  \label{dr:H}
H \colon \bigcup_{t \,\in\, [0,1]} \{ t \} \times \phi_t (U) \,\ra\,  \RR
\end{equation} 
generating the path $\phi_t$,
i.e., $\phi_t$ is the solution of the Hamiltonian system
\begin{equation}  \label{ee:Ham}
\left. \begin{array}{llll}
   \frac{d}{dt} \phi_t(z) &=& J \nabla H_t (\phi_t(z)),  
                                        &\; z \in U,\; t \in [0,1] , \\ [0.2em]
   \phi_0 (z)               &=& z,       &\; z \in U .
        \end{array}
   \right\}
\end{equation}
Here, $J$ denotes the standard complex structure defined by
\[
\omega_0 (z,w) = \langle Jz, w \rangle, \quad z,w \in \RR^{2n} .
\]
The function $H(z,t) = H_t(z)$ is determined by the first equation in
\eqref{ee:Ham} up to a smooth function $h(t) \colon [0,1] \ra \RR$. 
Notice that $0 \in \phi_t(U)$ for all $t$.
We choose $h(t)$ such that 
\begin{equation}  \label{ea:H0}
H_t(0)=0 \quad \text{ for all } t \in [0,1] .    
\end{equation}

\m
\noindent 
{\bf Step 3. Intermezzo: End of the proof of Proposition \ref{pa2:ear}}

\s
\noindent
Before proceeding with the proof of Theorem \ref{ta2:ear} we shall prove 
\text{Proposition \ref{pa2:ear}}.
Fix a subset $A$ of $U$ whose closure $\overline{A}$ in $\RR^{2n}$ is
contained in $U$. Since $U$ is bounded, the set $\overline{A}$ is
compact, and so the set
\[
K = \bigcup_{t \,\in\, [0,1]} \{ t \} \times \phi_t (\overline{A}) 
            \,\subset\, [0,1] \times \RR^{2n}
\]
is also compact and hence bounded. 
We therefore find a bounded neighbourhood $V$ of $K$ which is open in
$[0,1] \times \RR^{2n}$ and is contained in the set
$\bigcup_{t \,\in\, [0,1]} \{ t \} \times \phi_t (U)$.
By Whitney's Theorem, there exists a smooth function 
$f$ on $[0,1] \times \RR^{2n}$ which is equal to $1$ on $K$ and
vanishes outside $V$. Since $V$ is bounded,
the function $fH \colon [0,1] \times \RR^{2n} \ra \RR$ has compact
support, and so the Hamiltonian system associated with $fH$ can be
solved for all $t \in [0,1]$. We define 
$\Phi_A$ to be the resulting time-$1$-map.
Then $\Phi_A$ is a globally defined symplectomorphism of $\RR^{2n}$ with
compact support, and $\Phi_A |_A = \ff |_A$.
The proof of \text{Proposition \ref{pa2:ear}} is thus complete.

\m
\noindent 
{\bf Step 4. End of the proof of Theorem \ref{ta2:ear}}

\s
\noindent
If the set $U$ is not bounded, the subset $A \subset U$ does not need to 
be relatively compact, and so there might be no cut off $fH$ of $H$
whose Hamiltonian flow exists for all $t \in [0,1]$.
We therefore need to extend the Hamiltonian $H$ more carefully.
We shall first verify that our assumption \eqref{ea2:2} on $\ff$ implies 
that $\nabla H$ is linearly bounded.
Since we do not know a direct way to extend a linearly bounded gradient 
field to a linearly bounded gradient field, we shall then pass to the function
\[
G(t,w) \,=\, \frac{H(t,w)}{g \left( \left| w \right| \right)}
\]
where $g \left( \left| w \right| \right) = \left| w \right|$ for $\left|
w \right|$ large.
Our assumption that $U$ is a Lip\-schitz domain will imply that $G$ is
Lipschitz continuous in $w$ and can hence be extended to a Lipschitz
continuous function $\widehat{G}$ on $[0,1] \times \RR^{2n}$.
After smoothing $\widehat{G}$ in $w$ to $\widetilde{G}$ we shall obtain
an extension $\widetilde{H} (t,w) = g \left( \left| w \right| \right)
\widetilde{G}(t,w)$ whose gradient is linearly bounded.

\begin{lemma}  \label{la:0}
Let $L>0$ be the constant guaranteed by assumption \eqref{ea2:2}.
\begin{itemize}
\item[(i)]
$\left| \phi_t(z) - \phi_t(z') \right| \,\ge\, 
L \left| z-z' \right|  
\quad
\text{ for all }\, t \in\;]0,1] \text{ and } z, z' \in U .$
\item[(ii)]
$\left\| d \phi_t (z) \right\| \,\le\, \displaystyle\frac{1}{L}
\quad 
\text{ for all }\, t \in \;]0,1] \text{ and } z \in U .$
\end{itemize}
\end{lemma}

\proof
(i)
In view of definitions \eqref{de:phie} and \eqref{de:phi} we have
\begin{equation}  \label{da:idae}
\phi_t(z) \,=\, \frac{1}{\eta(t)} \ff \left( \eta(t) z \right)
\end{equation} 
for all $t \in\;]0,1]$ and $z \in U$.
Together with assumption \eqref{ea2:2} we find 
\begin{eqnarray*}
\left| \phi_t(z) - \phi_t(z') \right| 
    &=&   \frac{1}{\eta(t)} \left| \ff (\eta(t)z) - \ff (\eta(t)z') \right| \\
    &\ge& \frac{1}{\eta(t)} L \left| \eta(t)z - \eta(t)z' \right|  \\
    &=&   L \left| z-z' \right| .
\end{eqnarray*}
Assertion (i) thus follows.

\s
(ii)
We fix $t \in\;]0,1]$ and $z \in U$.
Following the proof of Proposition 2.20 in \cite{MS} we decompose the
linear symplectomorphism $d \phi_t(z)$ as 
\[
d \phi_t(z) \,=\, PQ
\]
where both $P$ and $Q$ are symplectic
and $P$ is symmetric and positive definite and $Q$ is orthogonal.
According to \cite[Lemma 2.18]{MS} the eigenvalues of $P$ are of the
form 
\[
0 \,<\, \ll_1 \,\le\, \ll_2 \,\le\, \dots \,\le\, \ll_n 
  \,\le\, \ll_n^{-1} \,\le\, \dots \,\le\, \ll_2^{-1} \,\le\, \ll_1^{-1} .
\]
Since $Q$ is orthogonal, we find
\begin{equation}  \label{ida:APl}
\left\| d \phi_t(z) \right\| \,=\, \left\| P \right\| \,=\, \ll_1^{-1} . 
\end{equation} 
Let $v_1$ be an eigenvector of $\ll_1$. In view of assertion (i) we have 
\[
\ll_1 \left| v_1 \right| \,=\, \left| d \phi_t(z) v_1 \right| \,\ge\, L \left|
v_1 \right|
\] 
and so $\ll_1^{-1} \le L^{-1}$.
This and the identity \eqref{ida:APl} yield
$\left\| d \phi_t(z) \right\| \le L^{-1}$, 
and so assertion (ii) follows.
\proofend

For $r >0$ we denote by $B(r)$ the closed $r$-ball around $0 \in \RR^{2n}$.
We choose $\ee >0$ so small that $B(\ee) \subset U$.
Finally, we abbreviate
\begin{equation}  \label{da:Ut}
U_t \,=\, U \cap B \left( \frac{\ee}{e}\, e^{1/t} \right), \quad\; t \in\;]0,1].
\end{equation} 

\begin{lemma}  \label{la:1}
(i) 
There exists a constant $C_1>0$ such that
\[
\left| \nabla H_t(w) \right| \,\le\, \frac{C_1}{t^2} \left| w \right|
\quad 
\text{ for all }\, t \in\;]0,1] \text{ and } w \in \phi_t(U) .
\]
(ii) There exists a constant $c_1>0$ such that
\[
\left| \nabla H_t(w) \right| \,\le\, 
\frac{c_1}{t^2} e^{-1/t} \left| w \right|  
\quad
\text{ for all }\, t \in\;]0,1] \text{ and } w \in \phi_t(U_t) .
\]
\end{lemma}
  
\proof
(i)
Fix $t \in\;]0,1]$ and $w = \phi_t(z)$.
Using the first line in \eqref{ee:Ham} and the definitions
\eqref{da:idae} and \eqref{de:eta} we compute
\begin{eqnarray}  
J \nabla H_t(w) &=& \frac{d}{dt} \phi_t (z) \notag \\ 
                &=& \frac{d}{dt} \left( \frac{1}{\eta(t)} \ff
                                     (\eta(t)z) \right) \notag \\
                &=& \frac{\eta'(t)}{\eta(t)} \left( - \frac{1}{\eta(t)} \ff
                       (\eta(t)z) + d \ff (\eta(t) z) z \right) \label{ide:diff} \\
                &=& \frac{2}{t^2} 
                       \,\big( -w + d \ff ( \eta(t) z ) z \big) . \label{ide:1}
\end{eqnarray} 
Lemma \ref{la:0}\,(ii) with $t=1$ yields
\begin{equation}  \label{re:dp}
\| d \ff (z) \| \,\le\, \frac{1}{L} \quad \text{ for all }\, z \in U 
\end{equation}
and the identity $\phi_t(0) = \frac{1}{\eta(t)} \ff(0) =0$ and Lemma
\ref{la:0}\,(i) with $z'=0$ yield
\begin{equation}  \label{ide:wLz}
|w| = \left| \phi_t(z) \right| \,\ge\, L |z| .
\end{equation}
In view of the identity \eqref{ide:1} and the estimates 
\eqref{re:dp} and \eqref{ide:wLz} we conclude 
\begin{eqnarray*} 
\left| \nabla H_t(w) \right| \,=\, \left| J \nabla H_t(w) \right| 
   &\le& \frac{2}{t^2} 
         \left( |w|+ \| d \ff ( \eta(t) z) \| |z| \right) \\ 
   &\le& \frac{2}{t^2} 
         \left( |w|+ \frac{1}{L^2} |w| \right) \\ 
   &=&   \frac{2}{t^2} \left( 1+\frac{1}{L^2} \right) |w| .
\end{eqnarray*} 
The constant $C_1 := 2\left(1+\frac{1}{L^2} \right)$ is as required.

\s
(ii)
By the choice of $\ee$, the smooth map
$\ff$ is $C^2$-bounded on $B(\ee)$, and so
Taylor's Theorem applied to $\ff \colon B(\ee) \ra \RR^{2n}$ and 
$d \ff \colon B(\ee) \ra \cl (\RR^{2n})$ 
guarantees constants $M_1$ and $M_2$ such that for
each $x \in B(\ee)$,
\[
\begin{array}{rcll}
\ff(x) &=& \ff (0) + d \ff (0)x + r(x)  &\text{ with }\, \left| r(x)
\right| \le M_1 \left| x \right|^2, \\ [0.2em]
d \ff(x) &=& d \ff (0) + R(x)           &\text{ with }\, \left\| R(x)
\right\| \le M_2 \left| x \right| ,
\end{array}
\]
where $\| R(x) \|$ denotes the operator norm of the linear operator
$R(x) \in \cl (\RR^{2n})$.
Since $\ff (0) = 0$ and $d \ff (0) = id_{\RR^{2n}}$ we conclude that
\[
\left| \ff (x) - d \ff (x) x \right| \,=\, \left| r(x) - R(x)x \right|
\,\le\, (M_1+M_2) \left| x \right|^2 \quad \text{if } \left| x \right| \le 
\ee
\]
and so, with $x = \eta(t)z$,
\begin{equation}  \label{e:M12}
\left| \frac{1}{\eta(t)} \ff \left( \eta(t)z \right) - d \ff \left(
\eta(t)z \right) z \right| 
\,\le\, (M_1+M_2) \eta(t) \left| z \right|^2 \quad \text{if } \eta(t)
\left| z \right| \le \ee . \;
\end{equation}  
Assume now $z \in U_t$.
In view of the definition \eqref{da:Ut} of $U_t$ we then have
\[
\eta(t) \left| z \right| \,\le\, e^2 e^{-2/t} \,\frac{\ee}{e}\, e^{1/t} \,=\,
e e^{-1/t} \ee \,\le\, \ee .
\]
Inserting the estimate \eqref{e:M12} into \eqref{ide:diff} and using
\eqref{ide:wLz} we conclude that 
\begin{eqnarray*}
\left| \nabla H_t(w) \right| 
    &\le& \frac{2}{t^2} (M_1 + M_2) \eta(t) \left|z\right|^2 \\
    &\le& \frac{2}{t^2} (M_1 + M_2) e e^{-1/t} \ee \left|z\right| \\
    &\le& \frac{2}{t^2} e^{-1/t} (M_1 + M_2) e \ee \frac{1}{L}  \left|w\right| .        
\end{eqnarray*}
The constant $c_1 := 2(M_1 + M_2) e\:\!\ee \tfrac{1}{L}$ is as required.
\proofend

\begin{lemma}  \label{la:2}
(i) 
There exists a constant $C_2>0$ such that
\[
\left| H_t(w) \right| \,\le\, \frac{C_2}{t^2} |w|^2 \quad
          \text{ for all }\, t \in\;]0,1] \text{ and } w \in \phi_t(U) .
\]
(ii) 
There exists a constant $c_2>0$ such that
\[
\left| H_t(w) \right| \,\le\, \frac{c_2}{t^2} e^{-1/t}|w|^2 \quad
          \text{ for all }\, t \in\;]0,1] \text{ and } w \in \phi_t(U_t) .
\]
\end{lemma}
  
\proof
(i)
Fix $t \in\;]0,1]$ and $w = \phi_t(z)$.
The smooth path 
\[
\gg \colon [0,1] \ra \phi_t(U), \qquad \gg(s) = \phi_t(sz)
\]
joins $0$ with $w$.
Since $H_t(0) =0$ we find that 
\begin{eqnarray}  \label{ide2:H}
H_t(w) &=& H_t(0) + 
             \int_0^1 \langle \nabla H_t(\gg(s)), \gg'(s) \rangle \,ds
                                                               \notag \\
       &=& \int_0^1 \langle \nabla H_t(\phi_t(sz)), d \phi_t(sz)z \rangle \,ds.
\end{eqnarray} 
The identity $\phi_t(0)=0$, the mean value theorem and Lemma \ref{la:0}\,(ii) yield
\begin{equation}  \label{este:Ms}
\left| \phi_t(sz) \right| \;=\; \left| \phi_t(sz) - \phi_t(0) \right|
           \;\le\; \frac{1}{L} s |z| .
\end{equation} 
Using the identity \eqref{ide2:H}, Lemma \ref{la:1}\,(i), 
Lemma \ref{la:0}\,(ii) 
and the estimates \eqref{este:Ms} and \eqref{ide:wLz} we
can estimate
\begin{eqnarray*} 
\left| H_t(w) \right| 
   &\le& \int_0^1 \left| \nabla H_t(\phi_t(sz)) \right| \left| d \phi_t(sz) z
                      \right| \,ds \\ 
   &\le& \frac{C_1}{t^2} \frac{1}{L} \,|z| \int_0^1 \left| \phi_t(sz) \right|
                            \, ds \\
   &\le& \frac{C_1}{t^2} \frac{1}{L^2} \,|z|^2\, \frac{1}{2} \\
   &\le& \frac{1}{2} \,C_1 \frac{1}{L^4} \frac{1}{t^2} \,|w|^2.
\end{eqnarray*} 
The constant $C_2 := \tfrac{1}{2} \;\!C_1 \frac{1}{L^4}$ is
as required.

\s
(ii)
Assume now $z \in U_t$.
Using Lemma \ref{la:1}\,(ii) and estimating as above we obtain
\[
\left|H_t(w)\right| \,\le\, \frac{1}{2} \,c_1 \frac{1}{L^4} 
             \frac{1}{t^2} \,e^{-1/t} \left| w \right|^2 .
\]
The constant $c_2 := \tfrac{1}{2}\, c_1 \frac{1}{L^4}$ is as required.
\proofend

Choose a smooth function $g \colon [0,\infty[ \; \ra [1,\infty[$ such
that
\begin{equation}  \label{da:gr}
g (r) \,=\, 
\left\{ \begin{array}{ll}
       1           & \text{ if }\; r \le \frac{1}{2} ,  \\ [0.2em]
       r           & \text{ if }\; r \ge 2
        \end{array}
   \right. 
\end{equation}
and $0 \le g'(r) \le 1$ for all $r$.

\begin{figure}[h] 
 \begin{center}
  \psfrag{1}{$1$}
  \psfrag{2}{$2$}
  \psfrag{r}{$r$}
  \psfrag{g}{$g(r)$}
  \leavevmode\epsfbox{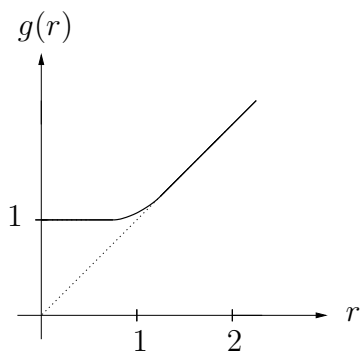}
 \end{center}
 \caption{The function $g(r)$.} \label{figurea2.fig}
\end{figure}
%
%magnification: 50
%

\ni
We define the smooth function $G \colon \bigcup_{t \,\in\, [0,1]} \{t\}
\times \phi_t(U) \ra \RR$ by
\begin{eqnarray}  \label{d:G}
G(t,w) \,\equiv\, G_t(w) \,:=\, \frac{H_t(w)}{g(|w|)} .
\end{eqnarray}
\begin{lemma}  \label{la:3}
(i)
There exists a constant $C_3>0$ such that
\[
\left| \nabla G_t(w) \right| \,\le\, \frac{C_3}{t^2} \quad
          \text{ for all }\, t \in\;]0,1] \text{ and } w \in \phi_t(U) .
\]
(ii)
There exists a constant $c_3>0$ such that
\[
\left| \nabla G_t(w) \right| \,\le\, \frac{c_3}{t^2} e^{-1/t} \quad
          \text{ for all }\, t \in\;]0,1] \text{ and } w \in \phi_t(U_t) .
\]
\end{lemma}
  
\proof
(i)
We have that
\[
\nabla \left( \frac{1}{g(|w|)}  \right) \,=\, - \frac{g'(|w|)}{g(|w|)^2}
\, \frac{w}{|w|}
\]
and so
\[
\nabla G_t(w) \,=\, - \frac{g'(|w|)}{g(|w|)^2} \, \frac{w}{|w|} \, H_t(w) +
\frac{1}{g(|w|)} \, \nabla H_t (w) .
\]
Using $g'(t) \in [0,1]$, Lemma \ref{la:2}\,(i) and Lemma \ref{la:1}\,(i) and $|w|
\le g(|w|)$ we therefore find that
\begin{eqnarray*} 
\left| \nabla G_t(w) \right|
     &\le& \frac{1}{g(|w|)^2} \left| H_t(w) \right| 
             + \frac{1}{g(|w|)} \left| \nabla H_t(w) \right| \\
     &\le& \frac{1}{g(|w|)^2} \, \frac{C_2}{t^2} \, |w|^2 
             + \frac{1}{g(|w|)} \, \frac{C_1}{t^2} \, |w| \\
     &\le& (C_1+C_2) \frac{1}{t^2} .
\end{eqnarray*} 
The constant $C_3 := C_1+C_2$ is as required.

\s
(ii)
Assume now $z \in U_t$. 
Using Lemma \ref{la:2}\,(ii) and Lemma \ref{la:1}\,(ii) and estimating
as above we obtain
\[
\left| \nabla G_t(w) \right| \,\le\, (c_1+c_2) \frac{1}{t^2} e^{-1/t}.
\]
The constant $c_3 := c_1+c_2$ is as required.
\proofend

\begin{lemma}  \label{la:4}
(i)
There exists a constant $C_4>0$ such that
\[
\left| G_t(w)- G_t(w') \right| \,\le\, \frac{C_4}{t^2} \left| w-w'
\right|
\quad
\text{ for all }\, t \in\;]0,1] \text{ and } w,w' \in \phi_t(U) .
\]
(ii)
There exists a constant $c_4>0$ such that
\[
\left| G_t(w)- G_t(w') \right| \,\le\, \frac{c_4}{t^2} e^{-1/t} 
\left| w-w' \right|
\quad
\text{ for all }\, t \in\;]0,1] \text{ and } w,w' \in \phi_t(U_t) .
\]
\end{lemma}  

\proof
(i)
Fix $t \in\;]0,1]$ and $w = \phi_t (z)$, $w' = \phi_t (z')$,
and assume that $U$ is a Lipschitz domain with Lipschitz constant $\ll$.
We then find a smooth path $\gg \colon [0,1] \ra U$ such that $\gg (0) = 
z$, $\gg (1) = z'$ and such that
\begin{eqnarray}  \label{ea:Llamda}
\length (\gg) \,=\, \int_0^1 \left| \gg'(s) \right| ds 
                       \,\le\, 2 \ll \left| z'-z \right| .
\end{eqnarray} 
Using Lemma \ref{la:3}\,(i), Lemma \ref{la:0}\,(ii), the estimate
\eqref{ea:Llamda} and Lemma \ref{la:0}\,(i) we can estimate
\begin{eqnarray*}
\left| G_t(w')-G_t(w) \right| &=& 
   \left| \int_0^1 \langle \nabla G_t \left( \phi_t (\gg(s)) \right), d
   \phi_t (\gg (s)) \gg' (s) \rangle ds \right| \\
  &\le& \frac{C_3}{t^2}\, \frac{1}{L} \int_0^1 \left| \gg' (s) \right| ds  \\
  &\le& \frac{C_3}{t^2}\, \frac{1}{L} \,2 \ll \left| z'-z \right| \\
  &\le& \frac{C_3}{t^2}\, \frac{1}{L^2} \, 2 \ll \left| w'-w \right| .
\end{eqnarray*}
The constant $C_4 := 2 C_3 \frac{1}{L^2} \ll$ is as required.

\s
(ii)
Assume now $z, z' \in U_t$. 
Since $U$ is starlike, we can assume that the path $\gg$ chosen above
is contained in $U_t$.
Using Lemma \ref{la:3}\,(ii) and estimating as above we obtain
\[
\left| G_t(w')-G_t(w) \right| \,\le\, \frac{c_3}{t^2} e^{-1/t}
\frac{1}{L^2} \:\!2 \ll \left|w'-w \right| .
\]
The constant $c_4 := 2 c_3 \tfrac{1}{L^2} \ll$ is as required.
\proofend

Our next goal is to extend the function $G$ on $\bigcup_{t \,\in\,
[0,1]} \{ t \} \times \phi_t (U)$ to a continuous function
$\widehat{G}$ on $[0,1] \times \RR^{2n}$ having similar properties.
We shall need two auxiliary lemmata.

\begin{lemma}  \label{la:McShane}
{\rm (Mc\,Shane \cite{Mc}\footnote{I'm grateful to Urs Lang for
pointing out to me this reference.})}
Consider a subset $W$ of the metric space $(X, d)$ and a function $f
\colon W \ra \RR$ which is $\ll$-Lipschitz continuous.
Then the function $\overline{f} \colon X \ra \RR$ defined by
\[
\overline{f} (x) \,:=\, \inf \left\{ f(w) + \ll\, d (x,w) \mid w \in W \right\}
\]
is a $\ll$-Lipschitz continuous extension of $f$.
\end{lemma} 

\begin{lemma}  \label{la:lip}
Assume that $V$ is a subset of $\RR^{2n}$ which contains the origin and
that the function
$h \colon V \cup B(2r) \ra \RR$ is $\ll_V$-Lipschitz continuous on $V$
and $\ll_B$-Lipschitz continuous on $B(2r)$.
Then $h$ is $(2\ll_V + \ll_B$)-Lipschitz continuous on $V \cup B(r)$.
\end{lemma} 

\proof
Fix $w,w' \in V \cup B(r)$.
If $w,w' \in V$ or $w,w' \in B(2r)$ then by assumption
\[
\left| h(w) - h(w') \right| \,\le\, \max (\ll_V, \ll_B) \left| w-w'
\right| .
\]
So assume that $w \in V \setminus B(2r)$ and $w' \in B(r)$.
Then 
$\left| w' \right| \le r \le \tfrac{\left| w \right|}{2}$ and so
\[
\frac{\left| w \right|}{2} \,\le\, \left| w \right| - \left| w' \right| \,\le\, 
\left| w-w' \right| . 
\]
Since $0 \in V$ and $0 \in B(2r)$ we can now estimate
\begin{eqnarray*}
\left| h(w) - h(w') \right| &\le& \left| h(w) - h(0) \right| + 
                                  \left| h(w') - h(0) \right|  \\
       &\le& \ll_V \left| w \right| + \ll_B \left| w' \right| \\
       &\le& (2 \ll_V + \ll_B) \frac{\left| w \right|}{2}  \\
       &\le& (2 \ll_V + \ll_B) \left| w-w' \right|,
\end{eqnarray*}
and so $h$ is $(2 \ll_V + \ll_B)$-Lipschitz continuous on $V \cup B(r)$.
\proofend

\begin{lemma}  \label{la:ext}
There exists a continuous function
$\widehat{G} \colon [0,1] \times \RR^{2n} \ra \RR$ 
with the following properties.
\begin{itemize}
\item[(i)]
$\widehat{G} (t,w) = G (t,w) \quad
          \text{ for all }\, t \in [0,1] \text{ and } w \in \phi_t(U) .$
%\item[(ii)]
%$\widehat{G} (0,w) = 0 \quad
%          \text{ for all }\, w \in \RR^{2n} .$
\item[(ii)]
There exists a constant $C_5>0$ such that
\[
\left| \widehat{G}_t(w)- \widehat{G}_t(w') \right| \,\le\,
   \frac{C_5}{t^2} \left| w-w' \right|
\quad
\text{ for all }\, t \in\;]0,1] \text{ and } w,w' \in \RR^{2n} .
\]
\end{itemize}
\end{lemma}

\proof
We shall first construct a function $\widehat{G} \colon [0,1] \times
\RR^{2n} \ra \RR$ meeting assertions (i) and (ii), and shall then 
verify that $\widehat{G}$ is continuous.
 
\s
We set $\widehat{G} (0,w) =0$ for all $w \in \RR^{2n}$. 
Since 
$H \colon \bigcup_{t \,\in\, [0,1]} \{ t \} \times \phi_t (U) \ra
\RR$ is continuous, Lemma \ref{la:2}\,(ii) and the definition \eqref{da:Ut}
of $U_t$ imply that $H(0,w) =0$ for all $w \in U$.  
In view of definition \eqref{d:G} we therefore have $G(0,w)=0$ for all
$w \in U$, and so assertion (i) holds for $t=0$.

We now fix $t \in\;]0,1]$.
We define the number $R_t$ by
\begin{equation}  \label{da:Rt}
R_t \,=\, \frac{L}{2} \frac{\ee}{e}\, e^{1/t} .
\end{equation}
Fix $w = \phi_t (z) \in B(2R_t)$.
In view of the estimate \eqref{ide:wLz} and the definition \eqref{da:Rt} 
we have
\[
\left| z \right| \,\le\, \frac{\left| w \right|}{L} \,\le\,
\frac{2}{L} R_t \,=\, \frac{\ee}{e}\, e^{1/t} ,
\]
and so $z \in U_t$ in view of definition \eqref{da:Ut}.
Lemma \ref{la:4}\,(ii) therefore implies that the function $G_t$ is
$\tfrac{c_4}{t^2} e^{-1/t}$-Lipschitz continuous on $\phi_t(U) \cap B(2
R_t)$.
According to Lemma \ref{la:McShane} the function $\overline{G}_t \colon
B(2R_t) \ra \RR$ defined by
\begin{equation}  \label{da:Gbar}
\overline{G}_t (x) \,:=\, \inf \left\{ G_t(w) + \frac{c_4}{t^2}e^{-1/t} 
\left|x-w \right| \,\big|\, w \in \phi_t(U) \cap B(2R_t) \right\} \;\;
\end{equation}
is a $\tfrac{c_4}{t^2} e^{-1/t}$-Lipschitz extension of $G_t$ to
$B(2R_t)$.
In particular, the function 
$\overline{\overline{G}}_t \colon \phi_t(U) \cup B(2R_t) \ra \RR$, 
\begin{equation}  \label{da:Gbbar}
\overline{\overline{G}}_t (x) \,:=\, 
\left\{ \begin{array}{ll}
          G_t(x)             & \text{ if }\; x \in \phi_t(U) , \\ [0.1em]
          \overline{G}_t(x)  & \text{ if }\; x \in B (2R_t) ,
        \end{array}
   \right. 
\end{equation}
is well-defined.
According to Lemma \ref{la:4}\,(i), $\overline{\overline{G}}_t$ is
$\tfrac{C_4}{t^2}$-Lipschitz continuous on $\phi_t(U)$, and according to 
the above, $\overline{\overline{G}}_t$ is $\tfrac{c_4}{t^2}$-Lipschitz
continuous on $B(2R_t)$.
According to Lemma \ref{la:lip}, the restriction of
$\overline{\overline{G}}_t$
to $\phi_t(U) \cup B(R_t)$ is therefore $\tfrac{C_5}{t^2}$-Lipschitz
continuous where we abbreviated
\[
C_5 \,:=\, 2 C_4 + c_4 .
\]
Applying Lemma \ref{la:McShane} once more, we find that the function
$
\widehat{G}_t \colon \RR^{2n} \ra \RR
$
defined by
\begin{equation}  \label{da:Ghat}
\widehat{G}_t (x) \,:=\, \inf \left\{ \overline{\overline{G}}_t(w) + 
                         \frac{C_5}{t^2} \left|x-w \right| \,\big|\, 
                         w \in \phi_t(U) \cup B(R_t) \right\} \;\;
\end{equation}
is a $\tfrac{C_5}{t^2}$-Lipschitz extension of the restriction of
$\overline{\overline{G}}_t$ to $\phi_t(U) \cup B(R_t)$.
In particular, 
\[
\widehat{G}(t,w) \,=\, G(t,w) \quad
          \text{ for all }\, w \in \phi_t(U) .
\]
The function $\widehat{G} \colon ]0,1] \times \RR^{2n} \ra \RR$ thus
defined therefore meets assertion (i) for $t \in \;]0,1]$ and assertion
(ii).

\begin{figure}[h] 
 \begin{center}
  \psfrag{p}{$\phi_t(U)$}
  \leavevmode\epsfbox{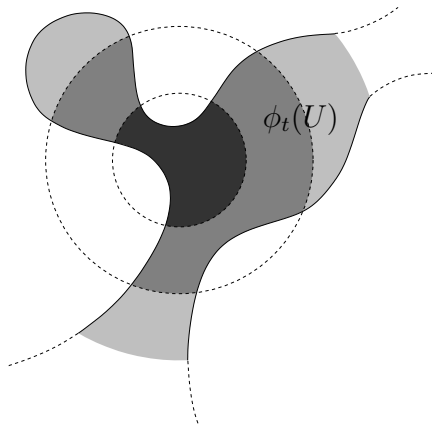}
 \end{center}
 \caption{The domain $\phi_t(U)$ and its intersections with 
             $B(R_t)$ and $B(2R_t)$.} \label{figurea3.fig}
\end{figure}
%
%magnification: 35
%

We are left with showing that the function 
$\widehat{G} \colon [0,1] \times \RR^{2n} \ra \RR$ constructed in the
previous two steps is continuous. 
The definitions \eqref{da:Rt}, \eqref{da:Gbar}, \eqref{da:Gbbar} and
\eqref{da:Ghat} show that the functions 
$\widehat{G}(\;\!\cdot\;\!, x) \colon ]0,1] \ra \RR$ and
$\widehat{G}(t,\:\!\cdot\:\!)  \colon \RR^{2n} \ra \RR$
are continuous.
This and the fact that the functions $\widehat{G}(t,\:\!\cdot\:\!)$ are
$\tfrac{C_5}{t^2}$-Lipschitz continuous imply that $\widehat{G} \colon
]0,1] \times \RR^{2n} \ra \RR$ is continuous.
In order to show that $\widehat{G}$ is also continuous at $(0,w)$ for
each $w \in \RR^{2n}$ we fix $w \in \RR^{2n}$.
We choose an open ball $B_w \subset \RR^{2n}$ centered at $w$. 
In view of definition \eqref{da:Rt} we have $R_t \ra \infty$ as $t \ra
0^+$.
We therefore find $t_0>0$ such that $B_w \subset B(R_t)$ for all $t \in
\;]0,t_0]$. We fix $t \in \;]0,t_0]$ and $w' \in B_w$.
Recalling the definition of $\widehat{G}(t,w') \,\equiv\, \widehat{G}_t(w')$
we see that
\begin{eqnarray*}
\widehat{G}_t (w') \;=\; \overline{\overline{G}}_t (w') &\!=\!&
\overline{G}_t(w')   \\
   &\!=\!& \inf \left\{ G_t(v) + \frac{c_4}{t^2} e^{-1/t} \left|w'-v \right|
\,\big|\, v \in \phi_t(U) \cap B(2R_t) \right\} .
\end{eqnarray*}
Since $0 = \phi_t(0) \in \phi_t(U) \cap B(2R_t)$ and $G_t(0) = H_t(0)
=0$ we conclude that
\begin{equation}  \label{ea:Gt0}
\widehat{G}_t(w') \,\le\, \frac{c_4}{t^2} e^{-1/t} \left| w' \right| .
\end{equation}
Moreover, we recall from the beginning of the proof of Lemma
\ref{la:ext} that $G_t$ is $\tfrac{c_4}{t^2} e^{-1/t}$-Lipschitz
continuous on $\phi_t(U) \cap B(2R_t)$.
This and $G_t(0) = 0$ yield 
\[
\left| G_t(v) \right| \,\le\, \frac{c_4}{t^2} e^{-1/t} \left| v \right| 
\quad \text{ for all }\, v \in \phi_t(U) \cap B(2R_t) .
\]
Therefore,
\begin{eqnarray*}
G_t(v) + \frac{c_4}{t^2} e^{-1/t} \left| w'-v \right| &\ge& 
            - \left| G_t(v) \right| + \frac{c_4}{t^2} e^{-1/t}
             \left|w'-v \right| \\
     &\ge& \frac{c_4}{t^2} e^{-1/t} ( - \left|v\right| +
                        \left|w'-v\right|) \\
     &\ge& - \frac{c_4}{t^2} e^{-1/t} \left| w' \right| 
\end{eqnarray*}
for all $v \in \phi_t(U) \cap B(2R_t)$.
We conclude that
\begin{equation}  \label{ea:Gt0'}
\widehat{G}_t(w') \,\ge\, - \frac{c_4}{t^2} e^{-1/t} \left| w' \right| .
\end{equation}
The estimates \eqref{ea:Gt0} and \eqref{ea:Gt0'}, which hold for all $t
\in \;]0,t_0]$ and $w' \in B_w$, now imply that
\[
\left| \widehat{G}_t(w') \right| \,\le\, \frac{c_4}{t^2} e^{-1/t} \left|
w' \right| 
\quad
   \text{ for all }\, t \in \;]0,t_0] \text{ and } w' \in B_w 
\]
and so $\widehat{G}$ is continuous at $(0,w)$. 
This completes the proof of Lemma \ref{la:ext}.
\proofend

Let now $A$ be a subset of $U$ whose closure in $\RR^{2n}$ is contained
in $U$.
Since also the origin is contained in $U$, we can assume that $A$ is
closed and $0 \in A$.
We abbreviate
\[
\ca \,:=\, \bigcup_{t \in [0,1]} \{t\} \times \phi_t(A) .
\]
The next step is to smoothen $\widehat{G}$ in the variable $w$ in such a 
way that the smoothened function $\widetilde{G}$ coincides with
$\widehat{G}$ on $\ca$.
We shall first construct a smooth function $G^*$ which approximates
$\widehat{G}$ very well and shall then obtain $\widetilde{G}$ by
interpolating between $\widehat{G}$ and $G^*$.

\s
Since $\RR^{2n}$ is a normal space, we find an open set $V$ in $\RR^{2n}$
such that
$A \subset V \subset \overline{V} \subset U$. Then
\begin{equation}  \label{ra:AVU}
\phi_t(A) \,\subset\, \phi_t(V) \,\subset\, \phi_t ( \overline{V} ) 
\,=\, \overline{\phi_t(V)} \,\subset\, \phi_t(U)  
\quad
          \text{ for all }\, t \in [0,1] .
\end{equation} 
We abbreviate
\[
\cv \,:=\, \bigcup_{t \in [0,1]} \{t\} \times \phi_t(V) .
\]
Since $\ca$ is closed and $\cv$ is open in $[0,1] \times \RR^{2n}$, 
we find a smooth function $f \colon [0,1] \times \RR^{2n} \ra [0,1]$ such
that
\begin{equation}  \label{ia:ft01}  
f |_\ca = 1 
\quad \text{ and } \quad
f |_{[0,1] \times \RR^{2n} \setminus \cv} = 0.
\end{equation}

We say that a continuous function $F \colon [0,1] \times \RR^{2n} \ra
\RR$ is 
{\it smooth in the variable $w \in \RR^{2n}$} if all derivatives $D^kF_t(w)$
of $F$ with respect to $w$ exist and are continuous on $[0,1] \times \RR^{2n}$.

\begin{lemma}  \label{la:star}
There exists a continuous function
$G^* \colon [0,1] \times \RR^{2n} \ra \RR$
which is smooth in the variable $w \in \RR^{2n}$ 
and has the following properties.
\begin{itemize}
\item[(i)]
$\left| \nabla f_t (w) \right| \left| G^*_t (w) - \widehat{G}_t (w)
\right| \,\le\, {\displaystyle\frac{C_5}{t^2}}  \quad
          \text{ for all }\, t \in\;]0,1] \text{ and } w \in \RR^{2n} .$
%\item[(ii)]
%$G^* (0,w) = 0 \quad
%          \text{ for all }\, w \in \RR^{2n} .$
\item[(ii)]
$
\left| \nabla G_t^*(w) \right| \le {\displaystyle\frac{2C_5}{t^2}}
\quad
\text{ for all }\, t \in\;]0,1] \text{ and } w \in \RR^{2n} .$
\end{itemize}
\end{lemma}

\proof
For each $l \in \NN$ we define the open subset $V_l$ of $[0,1] \times
\RR^{2n}$ by
\begin{equation}  \label{da:Vl}
V_l \,:=\, \left\{\:\! (t,w) \in [0,1] \times \RR^{2n} \mid \left| \nabla
f_t(w) \right| < l \:\!\right\} .
\end{equation}
Then there exists a smooth partition of unity 
$\left\{ \theta_i \right\}_{i \in \NN}$ 
on $[0,1] \times \RR^{2n}$ such that for each $i$ the support $\supp \theta_i$ 
is compact and contained in some $V_l$. We let $l_i$ be a number such
that $\supp \theta_i \subset V_{l_i}$.
Since $\left\{ \supp \theta_i \right\}$ form a locally finite covering of 
$[0,1] \times \RR^{2n}$, the set
\[
\Theta_i \,:=\, 
\left\{\:\! j \in \NN \mid \supp \theta_i \cap \supp \theta_j \neq
         \emptyset \:\!\right\}  
\]
is finite; let its cardinality be $m_i$.
We set
\begin{equation}  \label{da:Mi}
M_i \,:=\, \max \left\{\:\! m_j \mid j \in \Theta_i \:\!\right\} .
\end{equation}
Since the functions $\theta_i$ have compact support, the numbers 
\begin{equation}  \label{da:muit}
\mu_i \,:=\, \max \left\{ \left| \nabla \theta_i^t(w) \right| \,\big|\, (t,w) \in
[0,1] \times \RR^{2n} \right\} +1
\end{equation}
are finite.
We define positive numbers $r_i$ by
\begin{equation}  \label{da:rit}
r_i \,:=\, \frac{1}{l_i M_i \mu_i} .
\end{equation}
We next choose a smooth bump function $K \colon \RR^{2n} \ra [0,\infty[$ 
such that $\supp K \subset B(1)$ and $\int_{\RR^{2n}} K(v) \,dv =1$.
We abbreviate
\[
\kappa \,:=\, \max \left\{ \left| \nabla K (v) \right| \mid v \in
\RR^{2n} \right\} .
\]
For each $i$ we define a smooth function $K_i \colon \RR^{2n} 
\ra [0, \infty[$ by
\[
K_i (w) \,:=\, 
\frac{1}{r_i^{2n}} \,K \left( \frac{w}{r_i}\right) .
\]
Then $\supp K_i \subset B(r_i)$ and $\int_{\RR^{2n}} K_i (v) \,dv
=1$, and
\begin{equation}  \label{ea:gKti}
\left| \nabla K_i (w) \right| \,\le\, \frac{1}{r_i^{2n+1}}\, \kappa
\quad
          \text{ for all }\, w \in \RR^{2n} .
\end{equation}
Let $\widehat{G}$ be the function guaranteed by Lemma \ref{la:ext}.
For each $i$ we define the function $G_i^* \colon [0,1] \times \RR^{2n}
\ra \RR$ as the convolution
\begin{equation}  \label{da:Gistar}
G_i^* (t,w) \,:=\, \left( \widehat{G}_t * K_i \right) (w) \,\equiv\,
\int_{\RR^{2n}} \widehat{G}_t (v) K_i(w-v)\,dv .
\end{equation}
Since for each $t$ the function $\widehat{G}_t$ is continuous and since
$K_i$ is smooth, the function $w \mapsto G_i^*(t,w)$ is smooth and
\begin{equation}  \label{ia:Dk}
D^k G_i^*(t,w) \,=\, \int_{\RR^{2n}} \widehat{G}_t(v) D^k K_i(w-v) \,dv, 
\quad\; k=0,1,2, \dots
\end{equation} 
(see, e.g., \cite[Chapter 2, Theorem 2.3]{Hi}).
The function $\widehat{G}$ is continuous, and
$D^kK_i$ is continuous and has compact support and is thus uniformly continuous.
Formula \eqref{ia:Dk} therefore shows that $D^kG_i^*$ 
is continuous, $k=0,1,2, \dots$,
and so $G_i^*$ is continuous and smooth in $w$.
It follows that the function $G^* \colon [0,1] \times \RR^{2n} \ra \RR$
defined by
\begin{equation}  \label{da:Gstar}
G^* (t,w) \,:=\, \sum_i \theta_i (t,w) G_i^* (t,w)
\end{equation}
is continuous and smooth in $w$.
%Moreover, the identities $\widehat{G} (0,v) = 0$ holding for all 
%$v \in \RR^{2n}$ imply assertion (ii).
In order to prove assertions (i) and (ii) we fix $t \in\;]0,1]$ and abbreviate
\[
\theta_i(w) = \theta_i(t,w), \quad \:\!\!
\widehat{G}(w) = \widehat{G}(t,w), \quad \:\!\!
G_i^*(w) = G_i^*(t,w), \quad \:\!\!
G^*(w) = G^*(t,w).  
\]
%
%\noindent
{\it Proof of (i)}.\;
Using the definition \eqref{da:Gistar} of the function $G_j^*$ and the identity $\int_{\RR^{2n}} K_j
(v) \,dv =1$ we find
\begin{eqnarray*}
G_j^* (w) - \widehat{G}(w) 
&=& \int_{\RR^{2n}} \left( \widehat{G}(v) - \widehat{G}(w) \right) K_j (w-v)\, dv \\
&=& \int_{\RR^{2n}} \left( \widehat{G}(w-v) - \widehat{G}(w) \right) K_j
(v)\, dv 
\end{eqnarray*}
and so, together with Lemma \ref{la:ext}\,(ii),
\begin{eqnarray} \label{ea:2}
\left| G_j^* (w) - \widehat{G}(w) \right| 
&\le& 
\int_{\RR^{2n}} \left| \widehat{G}(w-v) - \widehat{G}(w) \right| K_j
(v)\, dv \notag\\
&\le& 
\int_{B(r_j)} \frac{C_5}{t^2} \left| v \right| K_j(v)\, dv \notag\\
&\le& 
\frac{C_5}{t^2} \, r_j \int_{B(r_j)}  K_j(v)\, dv \notag\\
&=&
 \frac{C_5}{t^2} \, r_j .
\end{eqnarray}
If $\nabla f_t(w) = 0$, assertion (i) is obvious. So assume $\left|
\nabla f_t (w) \right| > 0$.
Recall from the definitions \eqref{da:Mi} and \eqref{da:muit}
that $M_j \ge 1$ and $\mu_j \ge 1$.
This, the definition \eqref{da:rit} of $r_j$, the inclusion $\supp
\theta_j \subset V_{l_j}$ and the definition \eqref{da:Vl} of $V_{l_j}$ yield
\begin{equation}  \label{ea:rjL}
r_j\,=\, \frac{1}{l_j M_j \mu_j} \,\le\, \frac{1}{l_j} \,\le\,
\frac{1}{\left| \nabla f_t(w) \right|}
\quad
\text{ for all }\, w \in \supp \theta_j .
\end{equation}
The definition \eqref{da:Gstar} of $G^*$ and the estimates \eqref{ea:2} and
\eqref{ea:rjL} now yield
\begin{eqnarray*} 
\left| G^*(w) - \widehat{G}(w) \right| 
&=& 
\left| \sum_j \theta_j(w) \left( G_j^*(w) - \widehat{G}(w) \right)
\right| \\
&\le& 
\sum_j \theta_j(w) \frac{C_5}{t^2} \frac{1}{\left| \nabla f_t(w) \right|} \\
&=& 
\frac{C_5}{t^2} \frac{1}{\left| \nabla f_t(w) \right|}
\end{eqnarray*}
and so assertion (i) follows.

\s
\noindent
{\it Proof of (ii)}.\;
Using the definition \eqref{da:Gstar} of $G^*$ and the identities
$\sum_j \theta_j (w) = \sum_j \theta_j (w') = 1$ we compute that for all 
$w, w' \in \RR^{2n}$,
\begin{eqnarray} \label{ia:long}
G^* (w') - G^*(w) 
&=& 
\sum_j \theta_j (w') G_j^* (w') - \sum_j \theta_j (w) G_j^*(w) \notag\\
&=& 
\sum_j \left( \theta_j (w') - \theta_j (w) \right) \big( G_j^* (w') -
\widehat{G} (w') \big) \notag \\
& & \qquad \qquad \!\!  + \sum_j \theta_j (w) \left( G_j^* (w') - G_j^*(w) \right) . \qquad
\end{eqnarray}
Fix now $w$. We choose $i$ such that $\theta_i(w) >0$, and we choose an
open ball $B_w \subset \RR^{2n}$ centered at $w$ such that $B_w \subset
\supp \theta_i$.
Fix $w' \in B_w$.
In view of the mean value theorem and the definition \eqref{da:muit} of
$\mu_j$ we find that 
\begin{equation}  \label{ea:thmj}
\left| \theta_j (w') - \theta_j (w) \right| 
\,\le\, 
\max_{v \in B_w} \left| \nabla \theta_j (v) \right| \left|
w'-w \right| 
\,\le\,
\mu_j \left| w'-w \right|
\end{equation}
and the estimate \eqref{ea:2} with $w$ replaced by $w'$ yields
\begin{equation}  \label{ea:ggj}
\left| G_j^* (w') - \widehat{G} (w') \right| 
\,\le\, \frac{C_5}{t^2} r_j .
\end{equation}
The definition \eqref{da:Mi} of $M_j$ implies that
$M_j \ge m_i$ whenever $j \in \Theta_i$, and so
\begin{equation}  \label{ea:mi1}
\sum_{j \in \Theta_i} \frac{1}{M_j} \,\le\, \sum_{j \in \Theta_i}
\frac{1}{m_i} \,=\, 1 
\end{equation}
in view of the definition of $m_i$.
The definition \eqref{da:rit} of $r_j$ and the inequalities 
$l_j \ge 1$ and \eqref{ea:mi1} yield
\begin{equation}  \label{ea:mujr}
\sum_{j \in \Theta_i} \mu_j r_j \,=\, 
\sum_{j \in \Theta_i} \mu_j \frac{1}{l_jM_j\mu_j} \,\le\,
\sum_{j \in \Theta_i} \frac{1}{M_j} \,\le\,1 . 
\end{equation}
Since $w,w' \in B_w \subset \supp \theta_i$ we have
$\theta_j(w') - \theta_j(w) =0$ if $j \notin \Theta_i$.
This and the estimates \eqref{ea:thmj}, \eqref{ea:ggj} and
\eqref{ea:mujr} now show that
\begin{eqnarray} \label{ea:first}
\left| \sum_j \left( \theta_j (w') - \theta_j (w) \right) \left( G_j^* (w') -
     \widehat{G} (w') \right) \right| 
&\le& \sum_{j \in \Theta_i} \mu_j \left| w'-w \right| \frac{C_5}{t^2}\,
            r_j \notag\\ 
&\le& \frac{C_5}{t^2} \left| w'-w \right| .
\end{eqnarray}
Next, the definition \eqref{da:Gistar} of $G_j^*$ and the identity
$\int_{\RR^{2n}} K_j(v)\,dv =1$ yield
\begin{eqnarray*}
G_j^*(w') - G_j^*(w) 
&=& \int_{\RR^{2n}} \widehat{G}(v) \left( K_j (w'-v) - K_j(w-v) \right) dv \\
&=& \int_{\RR^{2n}} \left( \widehat{G}(w'-v) - \widehat{G} (w-v) \right) K_j (v)\,dv.
\end{eqnarray*}
Together with Lemma \ref{la:ext}\,(ii) we obtain
\[
\left| G_j^*(w') - G_j^*(w)  \right| \,\le\,
\frac{C_5}{t^2} \int_{\RR^{2n}} \left|w'-w \right| K_j(v)\,dv \,=\,
\frac{C_5}{t^2} \left|w'-w \right| 
\]
and so 
\begin{equation}  \label{ea:mid}
\left| \sum_j \theta_j(w) \left( G_j^*(w') - G_j^*(w) \right)  \right|
\,\le\,
\frac{C_5}{t^2} \left|w'-w \right| .
\end{equation}
The identity \eqref{ia:long} and the estimates \eqref{ea:first} and
\eqref{ea:mid} now imply
\[
\left| G^*(w') - G^*(w) \right| \,\le\, \frac{2C_5}{t^2} \left| w'-w
\right| .
\]
Since $w' \in B_w$ was arbitrary, we conclude that
\[
\left| \nabla G^* (w) \right| \,\le\, \frac{2C_5}{t^2}
\]
and so assertion (ii) follows.
The proof of Lemma \ref{la:star} is complete.
\proofend

\begin{lemma}  \label{la:smooth}
There exists a continuous function
$\widetilde{G} \colon [0,1] \times \RR^{2n} \ra \RR$ 
which is smooth in the variable $w \in \RR^{2n}$ 
and has the following properties.
\begin{itemize}
\item[(i)]
$
\widetilde{G} (t,w) = G (t,w) \quad
          \text{ for all }\, t \in [0,1] \text{ and } w \in \phi_t(A) .$
%\item[(ii)]
%$\widetilde{G} (t,w) = 0 \quad
%          \text{ for all }\, t \le 0 \text{ and } w \in \RR^{2n} .$
\item[(ii)]
There exists a constant $C_6>0$ such that
\[
\left| \nabla \widetilde{G}_t(w) \right| \,\le\,
   \frac{C_6}{t^2} 
\quad
\text{ for all }\, t \in\;]0,1] \text{ and } w \in \RR^{2n} .
\]
\end{itemize}
\end{lemma}

\proof
Let $f \colon [0,1] \times \RR^{2n} \ra [0,1]$ be the smooth function
chosen before Lemma \ref{la:star},
and let $\widehat{G}$ and $G^*$ be the continuous functions on $[0,1]
\times \RR^{2n}$ guaranteed by Lemma \ref{la:ext} and Lemma
\ref{la:star}.
We define a continuous function $\widetilde{G} \colon [0,1] \times
\RR^{2n} \ra \RR$ by
\[
\widetilde{G} (t,w) \,:=\,
f(t,w) \widehat{G} (t,w) + (1-f(t,w)) G^* (t,w).
\]
The inclusions \eqref{ra:AVU} and the identities \eqref{ia:ft01}, Lemma
\ref{la:ext}\,(i) and the fact that $G^*$ is smooth in $w$ imply that
$\widetilde{G}$ is smooth in $w$ and that assertion (i) holds true.
%Assertion (ii) follows from Lemma \ref{la:ext}\,(ii) and Lemma
%\ref{la:star}\,(ii). 
In order to verify assertion (ii) we fix $t \in\;]0,1]$. We first assume $w \in
\phi_t(U)$.
On $\phi_t(U)$ we have $\widehat{G}_t = G_t$, and so 
\[
\nabla \widetilde{G}_t(w) \,=\, 
\nabla f_t(w) \left( \widehat{G}_t(w)-G_t^*(w) \right)
+
f_t(w) \nabla G_t(w) + \left( 1-f_t(w) \right) \nabla G_t^* (w) .
\]
In view of Lemma \ref{la:star}\,(i), Lemma \ref{la:3}\,(i) and Lemma
\ref{la:star}\,(ii) we can therefore estimate
\begin{eqnarray*} 
\left| \nabla \widetilde{G}_t(w) \right| 
&\le&  
\left| \nabla f_t(w) \right| \left| G_t^*(w) - \widehat{G}_t(w) \right|
+
\left| \nabla G_t(w) \right| + \left| \nabla G_t^* (w)  \right| \\
&\le& 
\frac{C_5}{t^2} + \frac{C_3}{t^2} + \frac{2C_5}{t^2} .
\end{eqnarray*}
We next assume $w \in \RR^{2n} \setminus \overline{\phi_t(V)}$.
On $\RR^{2n} \setminus \overline{\phi_t(V)}$ we have $f_t \equiv 0$, and 
so 
\[
\left| \nabla \widetilde{G}_t(w) \right| \,=\, 
\left| \nabla G_t^*(w) \right| 
\,\le\, \frac{2C_5}{t^2} .
\]
Setting $C_6 := C_3 + 3 C_5$ assertion (ii) follows.
The proof of Lemma \ref{la:smooth} is complete.
\proofend

We are now in a position to define the desired extension $\widetilde{H}$
of $H$.
Let $g$ be the function chosen in \eqref{da:gr} and let $\widetilde{G}$
be the function guaranteed by Lemma \ref{la:smooth}.

\begin{lemma}  \label{la:Hti}
The continuous function
$\widetilde{H} \colon [0,1] \times \RR^{2n} \ra \RR$ defined by
\begin{equation}  \label{d:Ht}
\widetilde{H}(t,w) \equiv \widetilde{H}_t(w) \,:=\, g(\left|w\right|)
\widetilde{G}_t (w) .
\end{equation}
is smooth in the variable $w \in \RR^{2n}$ and has the following properties.
\begin{itemize}
\item[(i)]
$
\widetilde{H} (t,w) = H (t,w) \quad
          \text{ for all }\, t \in [0,1] \text{ and } w \in \phi_t(A) .$
%\item[(ii)]
%$\widetilde{H} (t,w) = 0 \quad
%          \text{ for all }\, t \le 0 \text{ and } w \in \RR^{2n} .$
\item[(ii)]
There exists a constant $C>0$ such that
\begin{equation}  \label{ea:Htti2}
\left| \nabla \widetilde{H}_t(w) \right| \,\le\, \frac{C}{t^2} \left(
\left|w\right|+1 \right)  \quad
          \text{ for all }\, t \in\;]0,1] \text{ and } w \in \RR^{2n} .
\end{equation}
\end{itemize}
\end{lemma}

\proof
Since $g \colon [0, \infty[ \;\ra [1,\infty[$ is smooth and
$\widetilde{G} \colon [0,1] \times \RR^{2n} \ra \RR$ is continuous, the
function $\widetilde{H}$ is indeed continuous, and since \text{$g(r) = 1$ if
$r \le \tfrac{1}{2}$} and $\widetilde{G}$ is smooth in $w$, the function
$\widetilde{H}$ is smooth in $w$.
Assertion (i) follows from the 
definition \eqref{d:Ht} of $\widetilde{H}$, from Lemma
\ref{la:smooth}\,(i) and from the definition \eqref{d:G} of $G$.
%and assertion (ii) follows from definition \eqref{d:Ht} and Lemma
%\ref{la:smooth}\,(ii).
In order to verify assertion (ii) we
fix $t \in\;]0,1]$ and $w \in \RR^{2n}$.
Using definition \eqref{d:Ht} we compute
\begin{equation}  \label{ea:grHti}
\nabla \widetilde{H}_t(w) \,=\, g'(\left| w \right|) \frac{w}{\left| w
\right|} \widetilde{G}_t(w) + g(\left| w \right|) \nabla
\widetilde{G}_t(w) .
\end{equation}
Since $0 \in A$ and $\phi_t(0) = 0$ we have, together with equation
\eqref{ea:H0},
\[
\widetilde{G}_t(0) \,=\, G_t(0) \,=\, H_t(0) \,=\, 0 .
\]
This, the mean value theorem and Lemma \ref{la:smooth}\,(ii) yield
\begin{equation}  \label{ea:G2x}
\left| \widetilde{G}_t(w) \right| \,\le\, \frac{C_6}{t^2} \left| w
\right| \quad \text{ and } \quad \left| \nabla \widetilde{G}_t(w)
\right| \,\le\, \frac{C_6}{t^2} .
\end{equation}
Using the identity \eqref{ea:grHti}, the estimates \eqref{ea:G2x} and
the estimates $\left| g'(r) \right| \le 1$ and $g(r) \le r+2$ holding
for all $r \ge 0$ we can estimate
\[
\left| \nabla \widetilde{H}_t (w) \right| \,\le\, \frac{C_6}{t^2}
\left| w \right| + ( \left| w \right| +2) \frac{C_6}{t^2} \,=\,
\frac{2C_6}{t^2} ( \left| w \right| +1 ) .
\]
Setting $C := 2C_6$ assertion (ii) follows.
The proof of Lemma \ref{la:Hti} is complete.
\proofend

Theorem \ref{ta2:ear} is a consequence of Lemma \ref{la:Hti}: 
The time-dependent vector field $\nabla \widetilde{H}_t(w)$ on $[0,1]
\times \RR^{2n}$ is continuous, and since it is smooth in $w$, it is
locally Lipschitz continuous in $w$.
This and assertion (ii) of Lemma \ref{la:Hti} imply that
the Hamiltonian system associated with $\widetilde{H}$ can be
solved for all $t \in [0,1]$. 
We define $\Phi_A$ to be the resulting time-$1$-map.
Since $\nabla \widetilde{H}_t(w)$ is continuous and smooth in $w$, the
map $\Phi_A$ is smooth (see \cite[Proposition 9.4]{A}), and so 
$\Phi_A$ is a globally defined symplectomorphism of $\RR^{2n}$.
Moreover, Lemma \ref{la:Hti}\,(i) shows that $\Phi_A |_A = \ff |_A$.
The proof of \text{Theorem \ref{ta2:ear}} is finally complete.
\proofend

\begin{remark}
{\rm
Proceeding as in Step 2 we obtain a {\it smooth}\, Hamiltonian 
\[
H_A \colon [0,1] \times \RR^{2n} \,\ra\, \RR
\]
which generates the symplectomorphism $\Phi_A$ and is such that $H_A
|_\ca = \widetilde{H} |_\ca$.
However, $H_A$ might not be $C^0$-close to $\widetilde{H}$, and $\nabla
H_A$ might not be linearly bounded.
}
\end{remark}


\begin{thebibliography}{99}


\bibitem{A} H.\ Amann. 
{\it Ordinary differential equations.}.
An introduction to nonlinear analysis.
de Gruyter Studies in Mathematics {\bf 13}. 
Walter de Gruyter \& Co., Berlin, 1990.

\bibitem{EH1} I.\ Ekeland and H.\ Hofer. Symplectic topology and
Hamiltonian dynamics. {\it Math. Z.} {\bf 200} (1990) 355--378.

\bibitem{Hi} M.\ Hirsch. {\it Differential topology}. 
Graduate Texts in Mathematics {\bf 33}.
Springer-Verlag, New York--Heidelberg, 1976. 

\bibitem{MS} D.\ Mc\;\!Duff and D.\ Salamon. 
{\it Introduction to Symplectic Topology}. 
Oxford Mathematical Monographs, Clarendon Press, 1995.

\bibitem{Mc} E.\ J.\ Mc\,Shane. 
Extension of range of functions.
{\it Bull. Amer. Math. Soc.} {\bf 40} (1934) 837--842.

\bibitem{Diss} F.\ Schlenk. 
{\it Embedding problems in symplectic geometry}.
Diss.\ ETH No.\ 14254.
Z\"urich 2001.


\end{thebibliography}
\enddocument